\documentclass[11pt]{article}

\usepackage{amsthm}
\usepackage{amsfonts}
\usepackage{amssymb}
\usepackage{verbatim}
\usepackage[all]{xy}
\usepackage{amsmath, amscd, amsthm}
\usepackage{longtable}
\usepackage{amscd}
\usepackage{mathrsfs}
\usepackage{graphicx}
\usepackage{latexsym}

\newtheorem{theorem}{Theorem}[section]
\newtheorem{lemma}[theorem]{Lemma}
\newtheorem{prop}[theorem]{Proposition}

\newtheorem{rem}[theorem]{Remark}

\usepackage{latexsym}

\DeclareMathOperator{\Wr}{Wr}
\DeclareMathOperator{\diag}{diag}
\begin{document}


\title{On Kaluzhnin - Krasner's embedding of groups}

\author{ A. Yu. Olshanskii \quad
\thanks{\ddag \quad The author was supported in part by the NSF grants DMS-1161294 and by the RFBR grant 15-01-05823.}
}
\maketitle

\begin{abstract}

In this note, we consider a 'thrifty' version of Kaluzhnin - Krasner's embedding
in wreath products and apply it to extensions by finite groups and to metabelian groups.
\end{abstract}


\medskip


\section{Introduction}

This note goes back to the pioneer paper of L. Kaluzhnin and M. Krasner \cite{KK},
where wreath products of groups were introduced and studied. Later many other group
theorists applied wreath products to construct various counter-examples and to
prove embedding theorems, and now wreath products are among the main tools of Group Theory. 
Here I pay attention to a feature of Kaluzhnin - Krasner's
works, that probably has not been used in  subsequent research papers. Herewith  I 
consider only standard wreath products of abstract groups (i.e., in terms of \cite{KK},  of permutation groups with regular actions). 

Let $A$ and $B$ be groups and $F$ a group of all functions $f: B\rightarrow A$ with
multiplication $(f_1f_2)(x)=f_1(x)f_2(x)$ for $x\in B$. The group $B$ acts on $F$ by shift
automorphisms: $(b\circ f) (x)= f(xb^{-1})$ for all $b, x\in B$, and the associated with
this action semidirect product $B\ltimes F$ is called the (complete) {\it wreath
product} of the groups $A$ and $B$, denoted by $A\Wr B$. Thus, every element
of $A\Wr B$ has a unique presentation as $bf$ ($b\in B, f\in F$) and the multiplication rule follows from the conjugation formula 
\begin{equation} \label{act}
(b^{-1}fb) (x) = f(xb^{-1})
\end{equation}
in $A\Wr B$ for any $b,x\in B$ and  $f\in F$.

Observe that any homomorphism $A\rightarrow \bar A$ induces the homomorphism $A\Wr B\rightarrow  \bar A\Wr B$ by the rule $bf\mapsto b\bar f$, where $\bar f\in \bar F$ is obtained by replacing the values of $f$ by their images in $\bar A$.

Given an arbitrary  group $G$ with a normal subgroup $A$, one has a canonical homomorphism
$\pi$ of $G$ onto the factor group $G/A = B$. Let $b\mapsto b^s$ be any
transversal $B\to G$, i.e. $\pi(b^s)=b$.  Then the Kaluzhnin - Krasner
monomorphism $\phi$ of the (abstract ) group $G$ into $A\Wr B$ is given by the
formula  (see \cite{KM}, \cite{N})
\begin{equation}\label{fun}
\phi(g) = \pi(g) f_g,\;\; where\;\; f_g(x) = (x\pi(g)^{-1})^s \;g \;(x^s)^{-1}
\end{equation}

Applying $\pi$ to $(x\pi(g)^{-1})^s \;g \;(x^s)^{-1}$, one obtains $1$, and so $f_g\in F$. To check that $\phi(g_1g_2)=\phi(g_1)\phi(g_2)$, one just exploits  the formulas (\ref{act},\ref{fun}). Finally,
$\phi$ is injective since obviously we have $\ker \phi \le A$, and by (\ref{fun}), 
$f_g(x)= x^s g (x^s)^{-1} \ne 1$ if $g\in A\backslash 1$.

The above-defined form of the Kaluzhnin - Krasner embedding $\phi$ is well known,
and, up to conjugation, the image $\phi(G)$ does not depend on the transversal $s$.
However, the original paper \cite{KK} suggested a stronger form of such an embedding.
Namely, assume now that a subgroup $C$ is normal in the normal subgroup $A$
of $G$, but $C$ does not contain nontrivial normal subgroups of $G$. Then consider
the product $\kappa$ of $\phi$ and the homomorphism $A\Wr B\rightarrow  \bar A \Wr B$
induces by the canonical homomorphism $A\rightarrow \bar A= A/C$ with $a\mapsto \bar a$ for $a\in A$. Thus
\begin{equation}\label{bar}
\kappa(g) = \pi(g) \overline{f}_g,\;\;   where\;\;
\overline{f}_g(x) = \overline{(x\pi(g)^{-1})^s \;g \;(x^s)^{-1}}
\end{equation}

If $g\in \ker\kappa$, then $\pi(g)=1$ and so $g\in A$.
Then formula (\ref{bar}) shows that the values $\overline{x^sg(x^s)^{-1}}$   are trivial
in $A/C$ for any $x\in B$, i.e. $x^sg(x^s)^{-1}\in C$. This implies
$hgh^{-1}\in C$ for every $h=ax^s\in G$, where $a\in A$, because $C$
is normal in $A$. It follows that $C$ contains the conjugacy class of $g$
in $G$, and so $g=1$ by the choice of the subgroup $C$.

Therefore $\kappa$ is also an embedding of $G$, and the following is
a slightly modified version of the old Kaluzhnin - Krasner statement:

\begin{prop} \label{kk} Let $G\rhd A \rhd C$ be a subnormal
series of a group $G$ with factors $\bar A=A/C$ and $B= G/A$, and let $C$ contain no non-trivial normal subgroups of $G$. Then there is an isomorphic embedding $\kappa$ of the
group $G$ in the wreath product $\bar A\Wr B= B \ltimes \bar F$, the embedding $\kappa$
can be defined by the formula (\ref{bar}),  so
$\kappa(A)\le \bar F$ and $\kappa(G)\bar F= \bar A\Wr B.$
\end{prop}

The group $\bar A$ can be much smaller than $A$, and making use of this below,
we 
\begin{itemize}
\item{apply the embedding $\kappa$ to wreath products with finite active group $B$,}
\item{embed finitely generated metabelian groups into $\bar A\Wr B$
with 'small' abelian $\bar A$ and $B$,  and}
\item{observe that $({\mathbb Z}/p{\mathbb Z})\Wr {\mathbb Z}$ and   ${\mathbb Z} \Wr  {\mathbb Z}$ contain $2^{\aleph_0}$ non-isomorphic
locally polycyclic subgroups.}
\end{itemize}

\section{Splittings of some group extensions}

The first application gives a characterization of wreath products with finite active group.

\begin{theorem} \label{odin} Assume that a normal subgroup $A$ of a group $G$ is a direct product
$\times_{i=1}^{n} H_i$, where $\{H_i\}_{1\le i\le n}$ is the set of all conjugate to
$H=H_1$ subgroups of $G$. Also assume that the normalizer $N_G(H)$ is equal to
$A$. Then $A$ has a semidirect compliment $B$
in $G$ and $G=H \Wr B$.
\end{theorem}
 
\proof 
Note that $n$ is the index of the normalizer of $H$ in $G$, therefore the group $B=G/A$ has order $n$. 

Define $C=C_1=\prod_{i\ne 1} H_i$. The normal in $A$
subgroup $C$ is conjugate in $G$ to any $C_j=\prod_{i\ne j} H_i$, as it follows from the assumptions. Therefore $\cap_{j=1}^n C_j=1$, and $C$ does not contain nontrivial
normal in $G$ subgroups.  

By Proposition \ref{kk}, we have the embedding $\kappa$ of $G$ in the wreath
product $(A/C)\Wr B$. If $g\in A$, then by formula (\ref{bar}),
$\kappa(g) =\bar f_g$, where $\bar f_g(x)= \overline{x^s g (x^s)^{-1}}.$ 
We may further assume that $x^s=1$ for $x=1$.
Then for $g\in H$, every value $x^s g (x^s)^{-1}$ of $f_g$ belongs to some $H_i\le C$ if $x\ne 1$ and $f_g(1)=g$. Therefore all the values of $\bar f_g$, except for one, are trivial, and
$\kappa (H)=H(1)$, where $H(1)$ is the subgroup of functions $\bar f$ supported by $\{1\}$ only. Since $H(1)\simeq A/C\simeq H_1=H$, the subgroup $H(1)$ can be identified with $H$.

Every conjugate of $H(1)$ in $\kappa(G)$ has the form $$(\pi(g')\bar f_{g'})H(1)(\pi(g')\bar f_{g'})^{-1}= bH(1)b^{-1},$$ where $b=\pi(g')$, because $H(1)$ is normal in $\bar F$. The subgroup $bH(1)b^{-1}$ consists of all functions $B\to H$ supported by $\{b^{-1}\}$. Here $b$ can be arbitrary element of the finite group $B$,
and therefore $\kappa(A)\ge \bar F$. The opposite inclusion $\kappa(A)\le \bar F$,
follows from Proposition \ref{kk}, that is $\kappa (A)=\bar F$.  Then
$A\Wr B= \kappa(G)\bar F =\kappa(G)$. Hence $\kappa$ is an isomorphism, and the theorem is proved. \endproof

The following example shows that one cannot remove the finiteness of the quotient
$B=G/A$ from the assumption of Theorem \ref{odin}.

\medskip

{\bf Example.} Let $G$ be a free metabelian group with two generators $a$ and $b$.
The commutator subgroup $A=[G,G]$ is the normal closure of one commutator $[a,b]$.
Since the subgroup $A$ is abelian, every subgroup conjugate to $H=\langle [a,b]\rangle$
is of the form $H_{ij}=c_{ij}^{-1}\langle [a,b]\rangle c_{ij}$, where $c_{ij} = a^ib^j$, $i,j\in\mathbb Z$. We explain below the known fact that $A$ is the direct product 
$\times_{i,j}H_{i,j}.$ It follows that $N_G(H) = A$. However $G$ does not split over $A$
because by  A. Shmelkin's \cite{Sh} result (also see \cite{N},
Theorem 42.56) two independent modulo $A$ elements must generate free metabelian
subgroup of $G$. 

To check that the elements $d_{ij}=c_{ij}^{-1}\langle [a,b]\rangle c_{ij}$ are linearly 
independent over $\mathbb Z$, one can apply the homomorphism $\mu$ of $G$
(in fact, a version of Magnus's embedding \cite{M, HM}) into the metabelian group of
upper triangle $2\times 2$ matrices over the group ring ${\mathbb Z}\langle x,y\rangle$ of a free abelian group of rank two given by the rule $$\mu (a)=\diag(x, 1), \;\; \mu(b)= \diag (y,1) + E_{12},$$
where  $E_{12}$ is the matrix with $1$ at position $(1,2)$ and zeros everywhere else. The matrix multiplication shows that $\mu(d_{ij})=I+x^{-i}y^{-j-1}(1-x^{-1})E_{12}$,
where the elements $x^{-i}y^{-j-1}(1-x^{-1})$ of ${\mathbb Z}\langle x,y\rangle$ ($i,j\in \mathbb Z$) are linearly independent over $\mathbb Z$.

\section{ Embeddings of metabelian groups}

There are finitely generated torsion free metabelian groups which are not embeddable
in $W={\mathbb Z} \Wr  B$ with finitely generated abelian $B$. For example, the derived subgroup $[G,G]$ of
the Baumslag - Solitar group $G=\langle a,b\mid b^{-1} b^{-1} ab= a^n\rangle$ is
isomorphic to the additive group of rationals whose denominators divide some powers 
of $n$, denote it by $D_n$; but for $|n|\ge 2$, $[W,W]$ has no nontrivial elements divisible by all powers of $n$. (Moreover, it is easy to construct $2$-generated metabelian groups $G$
 with $[G,G]$ containing an infinite direct powers of the group $D_n$ \cite{MO}.) However
the following is true.

\begin{theorem}\label{dva} Let $G$ be a finitely generated metabelian group with infinite
abelianization $B=G/[G,G]$. Then for some $n=n(G)\ge 1$,

(a)  $G$ embeds in the wreath product $ (D_n \times {\mathbb Z}/n {\mathbb Z})\Wr  B$.

(b) $G$ is isomorphic to a subgroup of $W=D_n \Wr  B$ if the derived subgroup $[G,G]$ is torsion free;

(c) $G$ is isomorphic to a subgroup of $W=({\mathbb Z}/n {\mathbb Z}) \Wr  B$, provided 
the derived subgroup $[G,G]$ is a torsion group;
\end{theorem}

The proof is based on the following 

\begin{lemma} \label{C} 
Let $G$ be a finitely generated metabelian group and $A$ an abelian normal subgroup of $G$. Then there is a subgroup $C$ in $A$ containing no non-trivial normal in $G$ subgroups, such that for some $n\ge 1$,

(a) the factor group $A/C$
is isomorphic to a subgroup of a finite direct power  of the group $D_n\times ({\mathbb Z}/n{\mathbb Z})$;

(b) $A/C$
is isomorphic to a subgroup of a finite direct power  of $D_n$ if $A$ is torsion free;

(c)  $A/C$
is isomorphic to a subgroup of a finite direct power  of ${\mathbb Z}/n{\mathbb Z}$ if $A$ is
a torsion group.

\end{lemma}

\proof {\bf (b)} We denote by $R$ the maximal normal in $G$ subgroup of finite (torsion-free) rank contained in $A$. Since $G$ satisfies the maximum condition for normal subgroups \cite{H}, this `radical' $R$ exists and contains all normal in $G$ subgroups of $A$ having finite rank.
Moreover, the maximal torsion subgroup $T/R$ of $L=A/R$ is trivial 
since the subgroup $T$ is normal in $G$ and its rank is equal to the rank of $R$.

We denote by ${\cal D}_n$ the class of all groups isomorphic to some subgroup of the countable direct power of the group $D_n$. 
By Ph.Hall's theorem of normal subgroups in metabelian groups (see \cite{H1}, we use the reformulation from \cite{MO}), we have $A\in {\cal D}_n$ for some $n$, and there exists a maximal linearly independent system $a_1, a_2,\dots $ in $A$ such that the torsion group $A/M$, where
$M=\langle a_1, a_2,\dots\rangle $, is a product of  Sylow $p$-subgroups for $p| n$ only.

Since 
$R$ has finite rank, it follows from the definition of of the class ${\cal D}_n$ that
the group $A\in {\cal D}_n$  has a subgroup $K$ such that the intersection $R\cap K$ is trivial and
the factor group $A/K$ is isomorphic to a subgroup of a finite direct power of $D_n$.

Let us enumerate non-trivial normal subgroups $N$ of $G$ of infinite rank, which are contained in $A$ and have torsion free factor groups $A/N$:
$N_1, N_2,\dots$ (This set is countable since $G$
satisfies the maximum condition for normal subgroups.) Using  this enumeration,
we will transform the basis of $M$ as follows.

Let $a_{01}=a_1, a_{02}=a_2,\dots $, and assume that the basis $(a_{i-1,1}, a _{i-1, 2},\dots)$
of $M$ is defined for some $i\ge 1$, and $a_{i-1,k}= a_k$ if $k$ is greater than  $m(i-1)$, where $m(0)=0$.
We will use additive notation for $A$.

Since the subgroup $N_i$ has infinite rank, it has a non-zero element 
$g_i=\sum_{j>m(i-1)}\lambda_j a_j$. Let $m(i)$ be the maximal subscript at non-zero coefficients $\lambda_j$ of this sum.
We may assume that the greatest common divisor of the coefficients $\lambda_{m(i-1)+1},\dots,  
\lambda_{m(i)} $ is $1$ since the group $A/N_i$ is torsion free.
Therefore there exists a basis $(a_{i,m(i-1)+1},\dots, a_{i,m(i)})$ of the free
abelian subgroup $\langle a_{m(i-1)+1},\dots, a_{m(i)}\rangle $ with $a_{i, m(i-1)+1}=g_i$. The other elements of the $(i-1)$-th basis of $M$ are left unchanged, i.e. $a_{i,k}=a_{i-1,k}$ if $k\le m(i-1)$ or $k>m(i)$.

Now we define $e_k$ = $a_{i,k}$ if $m(i-1)<k\le (m(i))$ for some $i$ and obtain a new
basis of $(e_1, e_2,\dots)$  of $M$ because we see from the construction that $\langle e_1,\dots, e_{m(i)}\rangle = \langle a_1, \dots, a_{m(i)}\rangle$ for every $m(i)>0$. Thus every
subgroup $N_i$ contains an element $g_i=a_{i,m(i-1)+1}=e_{m(i-1)+1}$ from this new basis.

Since the group $A$ is torsion free, every element of $A$ has a unique finite presentation
of the form $\sum_i\lambda_i e_i$ with  rational coefficients (although not every such a rational
combination belongs to $A$). Hence the subgroup $H$ of $A$ given by the equation $\sum_i\lambda_i=0$ is well defined, and the factor group $A/H$ is torsion free and has rank $1$. 
Note that $(M+H) /H\simeq M/(H \cap M)$ is infinite cyclic and $A/(M+H)$ is the factor group of the torsion group $A/M$. Therefore  the group $A/H$ is isomorphic to a subgroup of $D_n$.

It follows from the definition of $H$ that for every $i$, the subgroup $H$ does not contain any non-zero multiple $mg_i$ of the element $g_i=e_{m(i-1)+1}$.

Finally we set $C=H\cap K$. Then $A/C$ is embeddable in a finite direct power of $D_n$ 
since both $A/H$ and $A/K$ are. Assume now that $C$ contains
a nontrivial normal in $G$ subgroup $L$. Then $L$ has to have infinite rank since
$C\cap R\le K\cap R=0$. Let $T/L$ be the torsion part of $A/L$. Then $T$ is a normal
in $G$ subgroup with  torsion free factor group $A/T$. Hence $T=N_i$ for some $i$. Therefore
$L$ and $C$ must contain a non-zero multiple of the element $g_i\in N_i$. But $H$ does not contain such elements, and the statement (b) is  proved by contradiction.

{\bf (c)} The subgroup $A$ is a direct sum of its Sylow $p$-subgroups $A_p$. For every prime $p$, the elements $x$ of $A_p$ with  $px=0$ form a normal in $G$ subgroup $A(p)$.
It follows from the maximum condition for normal subgroups of $G$ that there are only finitely
many primes $p$ with nonzero  $A(p)$. For the same reason, $A$ has a finite exponent $n$.
Arguing as in the proof of (b), but replacing $A$ by $A(p)$ and taking the maximal 
set  $(a_1,a_2,\dots)$ in $A(p)$ as there (but linearly independent over ${\mathbb Z}/p{\mathbb Z}$),
we obtain a subgroup (and subspace) $C_p$ of finite codimension in $A(p)$, which does not contain any
nonzero normal in $G$ subgroup. 

Now consider a maximal  in $A_p$ subgroup $E_p$ such that $E_p\cap A(p)= C_p$.
Then every element $x+E_p$ of order $p$ from $A_p/E_p$ must belong to the canonical
image $(A(p)+E_p)/E_p$ of the subgroup $A(p)$ in $A_p/E_p$ since otherwise $(\langle x+E_p\rangle/E_p) \cap ((A(p)+E_p)/E_p)$ is trivial, and so \\$\langle  x  +E_p\rangle\cap A(p) \le E_p$, contrary to the maximality of $E_p$.
Since the  subgroup $(A(p)+E_p)/E_p\simeq A(p)/C_p$ is finite, we have finitely many
elements of order $p$ in the $p$-group $A_p/E_p$ of finite exponent dividing $n$. Hence $A_p/E_p$
is a finite $p$-group. 

If $N$ is a non-trivial $p$-subgroup normal in $G$  and $N\le E_p$, then 
$N\cap A(p) \le E_p\cap A(p) =C_p$, where $N\cap A(p)$ is non-trivial and normal in $G$,
contrary to the choice of $C_p$.
Thus $E_p$ contains no such subgroups $N$.

Since $A_p$ is a direct summand of $A$, for every $E_p$, one can find a subgroup $F_p\le A$ with $A_p\cap F_p = E_p$
and $A/F_p\simeq A_p/E_p$. If a normal in $G$ subgroup $N$ with nontrivial $p$-torsion 
were contained in $F_p$, then $A_p\cap N\le  A_p\cap F_p=E_p$, which would provide a contradiction.

Now the intersection $C= \cap_{p|n} F_p$ contains no nonzero  subgroups of $A$,
which are normal in $G$. Since every $A/F_p$ is a finite group of exponent dividing $n$,
the group $A/C$ is embeddable in a finite direct power of the group ${\mathbb Z}/n{\mathbb Z}$, as desired.

{\bf(a)} Let $T$ be the torsion subgroup of $A$. By the statement (b) applied to $G/T$, we have a subgroup $C'$ containing
$T$ but containing no bigger normal in $G$ subgroups, and $A/C'$ is embeddable in
a finite direct power of some $D_n$. Note that $C'$ contains no nontrivial torsion free
subgroup $N$ normal in $G$ since $N+T > T$.

Since $T$ has a finite exponent $m$, it has a torsion free
direct compliment $K$ in $A$ by Kulikov's theorem \cite{F,R}. The intersection $S$ of
the subgroups conjugated to $K$ in $G$ is torsion free, normal in $G$, and the exponent of $A/S$
is equal to the exponent of $A/K\simeq T$, i.e. it is $m$.

The statement (c) applied to $G/S$ provides us with a subgroup $C''$ 
containing
$S$ but containing no bigger normal in $G$ subgroups, with $A/C''$ embeddable in
a finite direct power of ${\mathbb Z}/m{\mathbb Z}$. Note that $C''$ contains no nontrivial torsion subgroup $N$ normal in $G$ since $N+S > S$.

It follows that $C=C'\cap C''$ contains no nontrivial normal in $G$ subgroups
and $A/C$ is embeddable in a finite direct power of $D_n\times {\mathbb Z}/m{\mathbb Z}$.
Since both $m$ and $n$ can be replaced by their
common multiple, the lemma is proved.
\endproof

\begin{rem} One cannot generalize Lemma \ref{C} to the slightly larger class of
central-by-metabelian groups. Indeed, every countable abelian group $A$ embeds as
a central subgroup in some finitely generated central-by-metabelian group $G$ \cite{H},
and so every subgroup $C$ of $A$ becomes normal in $G$. 

\end{rem}

{\bf Proof of Theorem \ref{dva}.}
{\bf (b)} By Lemma \ref{C} (b) and Proposition \ref{kk}, for some $n, m\ge 1$, the group $G$ embeds in a wreath product $W= D\Wr B$, where $D$ is a direct product of $m$ copies of the group $D_n$. Since $B$ is infinite, it follows from the Fundamental Theorem of finitely generated abelian groups, that $B$ contains a subgroup of index $m$ isomorphic to $B$. In other words, $B$ is a subgroup of index $m$
in a group $B_0$ isomorphic to $B$. The group
$W_0=D_n\Wr  B_0=B_0\ltimes F$ , where $F$ is the subgroup of functions $B_0\to D_n$,
contains the subgroup $W_1=B\ltimes F$, and it remains to show that
$W_1$ is isomorphic with $W$. This isomorphism is identical on $B$ and maps every function
$f_0: B_0\to D_n$ to the function $ f: B\to D$  given by the rule $f(b) = (f(t_1b),\dots, f(t_mb))\in D$,
where $\{t_1,\dots, t_m\}$ is a transversal to the subgroup $B$ in $B_0$.

{\bf (c,a)} One should argue as in the proof of (b) but with reference to the items (c) and (a)
of Lemma \ref{C} and  with the group $D_n$ replaced 
by ${\mathbb Z}/n{\mathbb Z}$ and $D_n\times {\mathbb Z}/n{\mathbb Z}$, respectively.

\section{Subgroups of $({\mathbb Z}/p{\mathbb Z})\;Wr\; {\mathbb Z}$ }

Let us fix a prime number $p$. For every prime $q\ne p$, there is a finite-dimensional
faithful, irreducible presentation of ${\mathbb Z}/q{\mathbb Z}$ over the Galois
field ${\mathbb F}_p$. Let $V_q$ be the corresponding representation module.
Each of these representations lifts to a representation of an infinite cyclic group $\langle b \rangle$, so the direct sum $V=\oplus_{q\ne p}V_q$ is a ${\mathbb F}_p\langle b \rangle$-module. The action of $\langle b \rangle$  defines the semidirect product $G=\langle b \rangle\ltimes V$.

\begin{lemma} \label{G} (1) Every finitely generated subgroup of $G$ is finite-by-cyclic,
i.e. $G$ is  locally finite-by-cyclic.

(2) $G$ contains $2^{\aleph_0}$ non-isomorphic subgroups.

(3) There is a subgroup $C$ in $V$ such that $V/C$ has order $p$ and
$C$ contains no nontrivial normal in $G$ subgroups.
\end{lemma}

\proof {\bf (1)} It is easy to see that any finite subset of $G$ is contained in
a subgroup $\langle b \rangle\ltimes U$, where $U$ is the direct sum of the 
subgroups $V_{q_i}$ over a finite subset of prime numbers $q_i$.
Since $U$ is finite, the statement (1) follows.

{\bf (2)} Denote by $H_S$ the subgroup $\langle b \rangle\ltimes V_S$, where
$S$ is a set of prime numbers $q\ne p$ and $V_S= \oplus_{q\in S} V_q$.
Since different $V_q$-s are irreducible and non-isomorphic ${\mathbb F}_p\langle b \rangle$-modules, every normal in $H_S$ finite
subgroup of $V_S$ is  $V_T= \oplus_{q\in T} V_q$ for a finite subset $T\subset S$.
The centralizer of $V_T$ in $H_S$ has index $\prod_{q\in T} q$. It follows that the
groups $H_{S_1}$ and $H_{S_2}$ are not isomorphic for $S_1\ne S_2$, which implies
the statement (2).

{\bf (3)} Let form an ${\mathbb F}_p$-basis $(e_1, e_2,\dots)$ as the union of
the bases of the subspaces $V_q$ and define the subspace $C$ by the equation
$\sum_i x_i=0$ in the coordinates. Then clearly $V/C$ has order $p$ and
$C$ contains none of the summands $V_q$. Since every normal in $G$ subgroup of $V$
is a submodule of the direct sum of some non-isomorphic irreducible ${\mathbb F}_p\langle b \rangle$-modules $V_q$, it
coincides with some $V_S$, and the statement is proved.
\endproof

The next theorem demonstrates that the structure of  subgroups of wreath products of
cyclic groups is quite rich.

\begin{theorem}\label{tri} The wreath product ${\mathbb Z}_p\Wr {\mathbb Z}$ contains
$2^{\aleph_0}$ non-isomorphic subgroups $H$, where each $H$ 
is locally finite-by-cyclic. 
\end{theorem} 

\proof The group $G$ from Lemma \ref{G} is embeddable in ${\mathbb Z}_p\Wr {\mathbb Z}$ 
by the property (3) of Lemma \ref{G} and Proposition \ref{kk}. Therefore Theorem \ref{tri} follows from
the properties (1) and (2) of Lemma \ref{G}.
\endproof

\begin{rem} Similar approach shows that  the wreath product ${\mathbb Z}\Wr {\mathbb Z}$ contains $2^{\aleph_0}$ countable locally polycyclic, non-isomorphic subgroups.  But now,
instead of $V_q$,
one should start with ${\mathbb Z}\langle b\rangle$-modules $V_i$, which are non-isomorphic, finite-dimensional and  irreducible over $\mathbb Q$. 
\end{rem}


\vskip11mm

Alexander Yu. Olshanskii:\newline
Department of Mathematics\newline
Vanderbilt University\newline
Nashville, TN 37240, USA.\newline
E-mail: alexander.olshanskiy@vanderbilt.edu

\end{document}